\documentclass[aop,MSNbibl,seceqn,dvips]{arximspdf}
\usepackage{mathrsfs}

\doi{10.1214/13-AOP900} 
\volume{42}
\issue{5}
\pubyear{2014}
\firstpage{1885}
\lastpage{1910}

\makeatletter
\newcommand{\rrVert}{\Vert}
\newcommand{\rrvert}{\vert}
\newcommand{\llVert}{\Vert}
\newcommand{\llvert}{\vert}
\newtheorem{theorem}{Theorem}[section]
\newtheorem{lemma}[theorem]{Lemma}
\newproclaim{remark}[theorem]{Remark}
\newproclaim{definition}[theorem]{Definition}
\newtheorem{proposition}[theorem]{Proposition}
\makeatother

\begin{document}
\begin{frontmatter}
\title{Densities for SDEs driven by degenerate\\ $\alpha$-stable processes}
\runtitle{Densities for degenerate SDEs with jumps}

\begin{aug}
\author[A]{\fnms{Xicheng} \snm{Zhang}\corref{}\ead[label=e1]{XichengZhang@gmail.com}\thanksref{t1}}
\runauthor{X. Zhang}
\affiliation{Wuhan University}
\address[A]{School of Mathematics and Statistics\\
Wuhan University\\
430072, Hubei\\
P.R. China\\
\printead{e1}} 
\end{aug}
\thankstext{t1}{Supported by NNSFs of China (Nos. 11271294, 11325105) and
Program for New Century Excellent Talents in University (NCET-10-0654).}

\received{\smonth{7} \syear{2012}}
\revised{\smonth{11} \syear{2013}}

%
\begin{abstract}
In this work, by using the Malliavin calculus, under H\"ormander's
condition, we prove the existence of
distributional densities for the solutions of stochastic differential equations
driven by degenerate subordinated Brownian motions. Moreover, in a
special degenerate case,
we also obtain the smoothness of the density. In particular,
we obtain the existence of smooth heat kernels for the following
fractional kinetic Fokker--Planck (nonlocal) operator:
\[
\mathscr{L}^{(\alpha)}_b:=\Delta^{\alpha/2}_\mathrm
{v}+\mathrm {v}\cdot \nabla_x +b(x,\mathrm{v})\cdot
\nabla_\mathrm{v},\qquad x,\mathrm{v}\in\mathbb{R}^d,
\]
where $\alpha\in(0,2)$ and $b\dvtx \mathbb{R}^d\times\mathbb{R}^d\to
\mathbb{R}^d$ is
smooth and
has bounded derivatives of all orders.
\end{abstract}



\pdfkeywords{60H07, 60H10, 35Q84, Malliavin calculus,
Hormander's condition, alpha-stable process,
distributional density, SDE}
\end{frontmatter}

\section{Introduction and main results}\label{sec1}
Consider the following stochastic differential equation (abbreviated as
SDE) in $\mathbb{R}^d$:
%
\begin{eqnarray}\label{Eq}
\mathrm{d}X_t&=&b(X_t)\,\mathrm{d}t+\sigma(X_t)
\,\mathrm{d}W_t+\int_{\mathbb{R}^d-\{0\}
}g(X_t,z)\widetilde N(\mathrm{d}t,\mathrm{d}z),
\nonumber\\[-10pt]\\[-10pt]
\eqntext{X_0=x,}
\end{eqnarray}
where $b\dvtx \mathbb{R}^d\to\mathbb{R}^d, \sigma\dvtx \mathbb{R}^d\to
\mathbb{R}^d\times\mathbb{R}^d$ and $g\dvtx \mathbb{R}^d\times\mathbb
{R}^d\to\mathbb{R}^d$ are smooth functions,
$(W_t)_{t\geq0}$ is a standard $d$-dimensional Brownian motion, and
$\widetilde N(\mathrm{d}t,\mathrm{d}z)$ is an independent compensated
Poisson random
measure on
$\mathbb{R}^d-\{0\}$ with intensity measure $\mathrm{d}t\nu(\mathrm{d}z)$.
Below, we always assume that $b,\sigma$ and $g$ have bounded
derivatives of all orders. Let us define the vector fields
\[
V_0:= \bigl(b_j-\tfrac{1}{2}\partial_l\sigma_{jk}\sigma_{lk} \bigr)
\partial_j\quad\mbox{and}\quad V_i:=\sigma_{ij}
\partial_j,\qquad i=1,\ldots, d,
\]
where we have used the convention: a repeated index in a product will
be summed automatically.
Set $\mathscr{V}_0:=\{V_1,\ldots, V_d\}$ and define recursively
\[
\mathscr{V}_k:= \bigl\{[V_0,V], [V_1,V],
\ldots,[V_d,V], V\in \mathscr{V}_{k-1} \bigr\},\qquad k\in
\mathbb{N},
\]
where $[V_i,V]:=V_iV-VV_i$ denotes the Lie bracket. It is well known
that when $g\equiv0$ (i.e., no jump part) and if $\bigcup_{k\in\mathbb
{N}}\mathscr{V}_k$ spans $\mathbb{R}^d$ at all points $x$ (called H\"
ormander's condition),
then the solution $X_t(x)$ of SDE (\ref{Eq}) admits a smooth density
$p_t(x,y)$, which was originally initiated by Malliavin \cite{Ma} (see
\cite{Nu} for a systematic introduction).
Moreover, by It\^o's formula, $p_t(x,y)$ satisfies the following
Fokker--Planck equation:
\[
\partial_tp_t(x,y)=\tfrac{1}{2}
\sigma_{ik}(x)\sigma _{jk}(x)\,\partial_{y_i}
\,\partial_{y_j}p_t(x,y)+\partial _{y_i}
\bigl(b_i(y)p_t(x,y)\bigr)
\]
with $p_0(x,y)=\delta_x(y)$.

Malliavin's probabilistic proof about H\"ormander's theorem is based on
the stochastic calculus of variations on the Wiener space invented by
himself \cite{Ma}.
Since then, there are many works devoted to extending the Malliavin
calculus to
the Poisson space case (see, e.g., \cite{Bi,Bi-Ja-Gr,Le,Pi,Ba-Cl},
etc.). In these works, the existence and smoothness of the
distributional densities for SDEs with jumps were obtained,
where various nondegeneracy conditions about $g(x,z)\nu(\mathrm{d}z)$
are imposed.
We particularly mention that Kusuoka in \cite{Ku} developed the
Malliavin calculus for subordinated Brownian motions,
and obtained the existence of smooth densities for SDEs driven by
nondegenerate subordinated Brownian motions. His argument will be
discussed later.

On the other hand, assuming that $\nu(\mathrm{d}z)=\mathrm
{d}z/|z|^{d+\alpha}$,
where $\alpha\in(0,2)$, and $g(x,z)$ satisfies
some boundedness and smoothness conditions, Takeuchi in \cite{Ta}, Corollary~1, proved that the solution $X_t(x)$ of SDE (\ref{Eq})
has a smooth density with respect to the Lebesgue measure
under some uniform H\"ormander's conditions. Notice that Takeuchi's
conditions allow pure-jump degenerate noises. In \cite{Ca}, Cass
obtained a similar result.
It is remarkable that recently, Kunita in \cite{Ku0} proved the
analytic property of distributional density to SDE (\ref{Eq}) under
weaker H\"ormander's conditions.
His proofs are based on the Malliavin calculus on the Wiener--Poisson
spaces developed in \cite{Is-Ku} and \cite{Ku1}. Moreover, an
estimate for discontinuous semimartingales
due to Komatsu and Takeuchi \cite{Ko-Ta} plays a crucial role in
Takeuchi and Kunita's proofs.
It is emphasized that all these results assume that $g$ is \textit{bounded} or the L\'evy measure $\nu$ has
\textit{finite} moments of all
orders. Thus,
the interesting $\alpha$-stable noise is ruled out.

In this work, we consider the following simple SDE:
%
\begin{equation}
\mathrm{d}X_t=b(X_t)\,\mathrm{d}t+A\,\mathrm{d}L_t,\qquad
X_0=x\in\mathbb {R}^d,\label{Eq0}
\end{equation}
where $A=(a_{ij})$ is a $d\times d$-matrix, and $(L_t)_{t\geq0}$ is a
rotationally invariant \mbox{$d$-}dimensional $\alpha$-stable process, that
is, its characteristic function
is given by
%
\begin{equation}
\mathbb{E}\mathrm{e}^{\mathrm{i}z\cdot L_t}=\mathrm {e}^{-t|z|^\alpha},\qquad \alpha
\in(0,2).\label{Eu4}
\end{equation}
We are interested in the problem that under what degenerate conditions
on $A$ together with $b$, $X_t(x)$ admits a smooth density with respect
to the Lebesgue measure.
Let us first look at the linear case of Ornstein--Uhlenbeck processes,
that is,
%
\begin{equation}
\mathrm{d}X_t=BX_t\,\mathrm{d}t+A\,\mathrm{d}L_t, \qquad X_0=x,\label{Eq1}
\end{equation}
where $B$ is a $d\times d$-matrix. The generator of this SDE is given
by $\mathcal{L}^{(\alpha)}_A+Bx\cdot\nabla$, where
the nonlocal operator $\mathcal{L}^{(\alpha)}_A$ is defined by
%
\begin{equation}
\mathcal{L}^{(\alpha)}_{A}f(x):=\mathrm{P.V.}\int
_{\mathbb
{R}^d}\bigl[f(x+A y)-f(x)\bigr]\frac{\mathrm{d}y}{|y|^{d+\alpha}},\label{EE1}
\end{equation}
where P.V. stands for the Cauchy principal value. Recently, Priola and
Zabczyk \cite{Pr-Za}
proved that $X_t$ has a smooth density under the following Kalman's
condition (see also \cite{Ku00}
for further discussions on this condition):
%
\begin{equation}
\operatorname{Rank}\bigl[A, BA,\ldots, B^{d-1}A\bigr]=d.\label{ER2}
\end{equation}
In fact, the solution of (\ref{Eq1}) is explicitly given by
\[
X_t=\mathrm{e}^{tB}x+\int^t_0
\mathrm{e}^{(t-s)B}A\,\mathrm{d}L_s=:\mathrm
{e}^{tB}x+Z_t.
\]
Using the approximation of step functions, by (\ref{Eu4}) it is easy
to see that
\[
\mathbb{E}\mathrm{e}^{\mathrm{i}z\cdot Z_t}=\mathbb{E}\exp \biggl\{ \mathrm{i}z\cdot\int
^t_0\mathrm{e}^{(t-s)B}A
\,\mathrm{d}L_s \biggr\} =\exp \biggl\{-\int^t_0\bigl|z^*
\mathrm{e}^{(t-s)B}A\bigr|^\alpha\,\mathrm{d}s \biggr\},
\]
where $*$ stands for the transpose of a column vector. Hence, for any
$m\in\mathbb{N}$,
\begin{eqnarray*}
\int_{\mathbb{R}^d}|z|^m\mathbb{E}\mathrm{e}^{\mathrm{i}z\cdot
Z_t}
\,\mathrm{d}z&=&\int_{\mathbb{R}^d}|z|^m\exp \biggl\{-\int
^t_0\bigl|z^*\mathrm{e}^{(t-s)B}A\bigr|^\alpha
\,\mathrm{d}s \biggr\}\,\mathrm{d}z
\\
&\leq&\int_{\mathbb{R}^d}|z|^m\exp \biggl
\{-|z|^\alpha\inf_{|a|=1}\int^t_0\bigl|a
\mathrm{e}^{sB}A\bigr|^\alpha\,\mathrm{d}s \biggr\} \,\mathrm{d}z.
\end{eqnarray*}
Here and below, ``$a$'' denotes a row vector in $\mathbb{R}^d$.
By (\ref{ER2}), one has
\[
\inf_{|a|=1}\int^t_0\bigl|a\mathrm{e}^{sB}A\bigr|^\alpha\,\mathrm{d}s>0
\]
and so,
\[
\int_{\mathbb{R}^d}|z|^m\mathbb{E}\mathrm{e}^{\mathrm{i}z\cdot
Z_t}
\,\mathrm{d}z<+\infty\qquad \forall m\in\mathbb{N}.
\]
Thus, $Z_t$ admits a smooth density by \cite{Sa}, Proposition~28.1,
and so does $X_t$.

We now turn to the nonlinear case. Before stating our main results, we
first recall some
notions about the subordinated Brownian motions. Let $(S_t)_{t\geq0}$
be a subordinator (an increasing one-dimensional
L\'evy process) on $\mathbb{R}_+$ with Laplace transform:
\[
\mathbb{E}\mathrm{e}^{-sS_t}=\exp \biggl\{t\int^\infty_0
\bigl(\mathrm {e}^{-s u}-1\bigr)\nu_S(\mathrm{d}u) \biggr\},
\]
where $\nu_S$ (called the L\'evy measure of $S_t$) satisfies $\nu_S(\{
0\})=0$ and
\[
\int^\infty_0(1\wedge u)\nu_S(
\mathrm{d}u)<+\infty.
\]
Below, we assume that $(S_t)_{t\geq0}$ is
independent of $(W_t)_{t\geq0}$ and
%
\begin{equation}
P \bigl\{\omega\dvtx  \exists t>0\mbox{ such that }S_t(\omega)=0 \bigr\}
=0,\label{CC1}
\end{equation}
which means that for almost all $\omega$, $t\mapsto S_t(\omega)$ is
\textit{strictly} increasing (see Lemma~\ref{EU1} below).
Notice that the Poisson process does not satisfy such an assumption,
but the $\alpha$-stable subordinator satisfies
this assumption (see \cite{Be}, p. 88, Theorem~11). Essentially,
condition (\ref{CC1}) is a nondegenerate assumption,
and says that the subordinator has infinitely many jumps on any
interval. In particular, the process defined~by
%
\begin{equation}
L_t:=W_{S_t},\qquad t\geq0,\label{EU2}
\end{equation}
is a L\'evy process (called subordinated Brownian motion) with
characteristic function:
\[
\mathbb{E}\mathrm{e}^{\mathrm{i} z\cdot L_t}=\exp \biggl\{t\int_{\mathbb{R}^d}
\bigl(\mathrm{e}^{\mathrm{i}z\cdot y} -1-\mathrm{i}z\cdot y1_{|y|\leq1}\bigr)
\nu_L(\mathrm{d}y) \biggr\},
\]
where $\nu_L$ is the L\'evy measure given by
%
\begin{equation}
\nu_L(\Gamma)=\int^\infty_0(2\pi
s)^{-d/2} \biggl(\int_\Gamma \mathrm{e}^{-|y|^2/2s}
\,\mathrm{d}y \biggr)\nu_S(\mathrm{d}s).\label{EW2}
\end{equation}
Obviously, $\nu_L$ is a symmetric measure.

The first aim of this paper is to prove the following existence result
of distributional density to SDE (\ref{Eq0})
under H\"ormander's condition as in \cite{Ta} and \cite{Ku0}.
%
%
\begin{theorem}\label{Th1}
Let $b\dvtx \mathbb{R}^d\to\mathbb{R}^d$ be a $C^\infty$-function with
bounded partial derivatives of first order.
For $x\in\mathbb{R}^d$, let $X_t(x)$ solve SDE (\ref{Eq0}) with
subordinated Brownian motion $L_t$.
Assume that for some $n=n(x)\in\mathbb{N}$,
{\renewcommand{\theequation}{$\mathscr{H}_n$}
%
\begin{equation}\label{eqHHN}
\operatorname{Rank}\bigl[A, B_1(x)A,
B_2(x)A,\ldots, B_n(x)A\bigr]=d,
\end{equation}}%
where $B_1(x):=(\nabla b)_{ij}(x)=(\partial_j b^i(x))_{ij}$, and for
$n\geq2$,
\setcounter{equation}{9}
%
\begin{equation}
B_n(x):=\bigl(b^i\,\partial_iB_{n-1}
\bigr) (x)-(\nabla b\cdot B_{n-1}) (x).
\end{equation}
Then the law of $X_t(x)$ is absolutely continuous with respect to the
Lebesgue measure.
In particular, the density $p_t(x,y)$ solves the following nonlocal
Fokker--Plack equation
in the weak or distributional sense:
%
\begin{equation}
\partial_tp_t(x,y)=\mathcal{L}_Ap_t(x,
\cdot) (y)+\partial _{y_i}\bigl(b_i(y)p_t(x,y)
\bigr)\label{ET8}
\end{equation}
with $p_0(x,y)=\delta_x(y)$, where
\[
\mathcal{L}_A f(y):=\mathrm{P.V.}\int_{\mathbb{R}^d}
\bigl[f(y+A z)-f(y)\bigr]\nu_L(\mathrm{d}z).
\]
\end{theorem}

%
\begin{remark}
If we assume that $L_t$ has finite moments of all orders, then this
result is contained in \cite{Ku0}, Theorem~5.1.
In fact, Kunita also obtained the smoothness of the density.
Nevertheless, our proof is simpler in this case. Notice that if $b(x)=B
x$, then condition (\ref{eqHHN}) reduces to (\ref{ER2}).
\end{remark}

For the smoothness of $p_t(x,y)$, we shall assume the following uniform
H\"ormander's condition:
{\renewcommand{\theequation}{$U\mathscr{H}_1$}
%
\begin{equation}\label{eqUH1}
\inf_{x\in\mathbb{R}^d}\inf_{|a|=1} \bigl(|aA|^2+\bigl|a
\nabla b(x)A\bigr|^2 \bigr)=:c_1>0
\end{equation}}\setcounter{equation}{11}%
and prove the following partial result.
%
%
\begin{theorem}\label{Th2}
Let $b\dvtx \mathbb{R}^d\to\mathbb{R}^d$ be a $C^\infty$-function with
bounded partial derivatives of all orders.
In addition to (\ref{eqUH1}), we assume that the L\'evy measure
$\nu_S$ satisfies for some $\theta\in(0,\frac{1}{2})$,
%
\begin{equation}
\lim_{\varepsilon\downarrow0}\frac{1}{\varepsilon^{1-2\theta
}}\int^\varepsilon_0
u\nu_S(\mathrm{d}u)=:c_\theta>0.\label{Con2}
\end{equation}
Then the density $p_t(x,y)$ is a smooth function on $(0,\infty)\times
\mathbb{R}^d\times\mathbb{R}^d$, and for each $t>0$,
\[
(x,y)\mapsto p_t(x,y)\in C^\infty_b\bigl(
\mathbb{R}^d\times\mathbb{R}^d\bigr).
\]
In particular, for all $(t,x,y)\in(0,\infty)\times\mathbb
{R}^d\times\mathbb{R}^d$,
%
\begin{equation}
\partial_t p_t(x,y)=\mathcal{L}_A
p_t(\cdot,y) (x)+b(x)\cdot\nabla _x
p_t(x,y).\label{Eq2}
\end{equation}
\end{theorem}
%
%
\begin{remark}
Condition (\ref{eqUH1}), compared with (\ref{eqHHN}), is
much stronger, and will be used to prove the
$L^p$-integrability of the inverse of the Mallavin covariance matrix
defined by (\ref{ETR1}) and (\ref{ETR2}) below,
where the key point is to prove a Norris' type lemma (see Lemma \ref
{Lemma2} below). We conjecture that a similar ($U\mathscr{H}_n$) as in
\cite{Ku0}
should imply the smoothness of $p_t(x,y)$. Nevertheless, the following
stochastic Hamilton system driven by a subordinated Brownian motion
satisfies (\ref{eqUH1}):
%
\begin{equation}
\label{Ham} \cases{ \mathrm{d}X_t=\nabla_y
H(X_t, Y_t)\,\mathrm{d}t, &\quad $X_0=x\in
\mathbb{R}^d$,
\vspace*{2pt}\cr
\mathrm{d}Y_t=-\nabla_xH(X_t,
Y_t)\,\mathrm{d}t+A\,\mathrm{d}L_t,&\quad
$Y_0=y\in\mathbb{R}^d$,}
\end{equation}
where $A$ is a $d\times d$-invertible matrix, and $H\dvtx \mathbb
{R}^d\times\mathbb{R}^d\to\mathbb{R}$ is a $C^2$-Hamiltonian
function so
that $y\mapsto H(x,y)$ is strictly convex or concave.
\end{remark}
%
%
\begin{remark}
Let $\nu_S(\mathrm{d}u)=u^{-(1+\alpha)}\,\mathrm{d}u$ be the L\'evy
measure of an
$\alpha$-stable subordinator. It is easy to see
that (\ref{Con2}) holds for $\theta=\alpha/2$.
\end{remark}

The argument for proving Theorems~\ref{Th1} and~\ref{Th2} is
different from Takeuchi and Kunita's works.
We shall follow Kusuoka's method \cite{Ku}. The advantage of which is
that it is not necessary to develop a \textit{new} Malliavin calculus for
jump processes,
and moreover, one can obtain some quantitive estimates about the
semigroup (see Theorem~\ref{Th4} below); while the drawback of which
is of course the loss of generality.
It is noticed that in \cite{Ku}, Kusuoka considered the SDE driven by
multiplicative noises. However, it seems that there is a gap in the
calculations about the
Malliavin covariance matrix (see \cite{Ku}, Theorem~3.3)
since the solution of SDE (\ref{Eq}) usually does not form a
stochastic diffeomorphism flow if there is no further restriction on
the jump size (cf. \cite{Pr}, p.~328).
This is also why we have to confine ourself to the additive noise.

Let us now describe the argument (see also \cite{Zh1}).
Let $(\mathbb{W},\mathbb{H},\mu_\mathbb{W})$ be the classical
Wiener space, that is, $\mathbb{W}$ is the space of all continuous
functions from $\mathbb{R}^+$ to~$\mathbb{R}^d$ with vanishing values at starting point $0$, $\mathbb
{H}\subset\mathbb{W}$ is the Cameron--Martin space consisting of all
absolutely continuous functions with square integrable derivatives, and
$\mu_\mathbb{W}$ is the Wiener measure so that
the coordinate process
\[
W_t(w):=w_t
\]
is a standard \mbox{$d$-}dimensional Brownian motion.

Let $\mathbb{S}$ be the space of all increasing, purely discontinuous
and c\`adl\`ag functions from $\mathbb{R}_+$ to $\mathbb{R}_+$ with
$\ell_0=0$, which is endowed with the Skorohod metric and
the probability measure $\mu_\mathbb{S}$ so that the coordinate process
\[
S_t(\ell):=\ell_t
\]
has the same law as the given subordinator. Consider the following
product probability space:
\[
(\Omega,\mathscr{F},P):=\bigl(\mathbb{W}\times\mathbb{S}, \mathscr {B}(
\mathbb{W})\times\mathscr{B}(\mathbb{S}), \mu_\mathbb {W}\times
\mu_\mathbb{S}\bigr)
\]
and define
\[
L_t(w,\ell):=w_{\ell_t}.
\]
Then $(L_t)_{t\geq0}$ has the same law as the given subordinated
Brownian motion. In particular, the solution
$X_t(x)$~of~SDE (\ref{Eq0}) can be regarded as a functional of $w$ and
$\ell$ and
%
\begin{equation}
\mathbb{E}f\bigl(X_t(x)\bigr)=\int_\mathbb{S} \int
_\mathbb{W}f\bigl(X_t(x,w_\ell )\bigr)
\mu_\mathbb{W}(\mathrm{d}w)\mu_\mathbb{S}(\mathrm{d}\ell
).\label{Ep6}
\end{equation}
The advantage of this viewpoint is that we can use the classical
Malliavin calculus to study
the Brownian functional $w\to X_t(x,w_\ell)$ (see \cite{Ku}).
Thus, in order to prove Theorem~\ref{Th1}, it is enough to prove that
for each $\ell\in\mathbb{S}$, the law of
$w\mapsto X_t(x,w_\ell)$ under $\mu_\mathbb{W}$ is absolutely
continuous with respect to the Lebesgue measure.
In order to prove Theorem 1.3, the key point is to prove the
$L^p$-integrability of the inverse of the Malliavin covariance matrix
so that we can use the integration by parts formula to derive some
gradient estimates (see Theorem~\ref{Th4} below),
which then implies the smoothness of the density by Sobolev's embedding theorem.

This paper is organized as follows: in Section~\ref{sec2}, we prove Theorem \ref
{Th1} by using the Malliavin calculus,
where the main point is to prove the invertibility of the Malliavin
covariance matrix
$(\Sigma^\ell_t)_{\ell=S}$ in (\ref{Ep4}) below. In Section~\ref{sec3}, we
prove Theorem~\ref{Th2}
by establishing a Norris' type lemma as in \cite{Ca}. In order to
overcome the nonintegrability of
$\alpha$-stable processes, we shall separately consider the small
jumps and the large jumps of the subordinator.
In particular, the asymptotic estimate of small times about the
semigroup plays a crucial role.


\section{Proof of Theorem \texorpdfstring{\protect\ref{Th1}}{1.1}}\label{sec2}
We need the following simple lemma about the density of the jump number
of the subordinator.
%
%
\begin{lemma}\label{EU1}
For $s>0$, set $\Delta\ell_s:=\ell_s-\ell_{s-}$ and
\[
\mathbb{S}_0:=\bigl\{\ell\in\mathbb{S}\dvtx  \{s\dvtx  \Delta
\ell_s>0\} \mbox{ is dense in $[0,\infty)$}\bigr\}.
\]
Under (\ref{CC1}), we have $\mu_\mathbb{S}(\mathbb{S}_0)=1$.
\end{lemma}
\begin{pf}
Let $\mathscr{I}$ be the total of all rational intervals in $[0,\infty
)$, that is,
\[
\mathscr{I}:=\bigl\{I=(a,b)\dvtx  0\leq a<b \mbox{ are rational numbers}\bigr\}.
\]
For $I\in\mathscr{I}$, let us write
\[
\mathbb{S}_I:=\bigl\{\ell\in\mathbb{S}\dvtx  I\subset\{s\dvtx  \Delta
\ell_s=0\} \bigr\}.
\]
It is easy to see that
\[
\mathbb{S}-\mathbb{S}_0=\bigcup_{I\in\mathscr{I}}
\mathbb{S}_I.
\]
Thus, for proving $\mu_\mathbb{S}(\mathbb{S}_0)=1$, it is enough to
prove that for each $I=\break (a,b)\in\mathscr{I}$,
\[
\mu_\mathbb{S}(\mathbb{S}_I)=\mu_\mathbb{S}\bigl(
\bigl\{\ell\in\mathbb {S}\dvtx  (a,b)\subset\{s\dvtx  \Delta\ell_s=0\}\bigr\}
\bigr)=0,
\]
which, by the stationarity of the subordinator, is equivalent to
%
\begin{equation}
\mu_\mathbb{S}\bigl(\bigl\{\ell\in\mathbb{S}\dvtx  (0,b-a)\subset\{s\dvtx  \Delta
\ell_s=0\}\bigr\}\bigr)=0.\label{Ep9}
\end{equation}
Since
\[
\bigl\{\ell\in\mathbb{S}\dvtx  (0,b-a)\subset\{s\dvtx  \Delta\ell_s=0\}\bigr
\}=\bigl\{\ell \in\mathbb{S}\dvtx  \ell_s=0,\ \forall s\in(0,b-a)\bigr\}
\]
by (\ref{CC1}), we obtain (\ref{Ep9}), and complete the proof.
\end{pf}

For a functional $F$ on $\mathbb{W}$, the Malliavin derivative of $F$
along the direction $h\in\mathbb{H}$ is defined as
%
\begin{equation}
D_hF(w):=\lim_{\varepsilon\to0}\frac{F(w+\varepsilon
h)-F(w)}{\varepsilon}\qquad\mbox{in }L^2(\mathbb{W},\mu_\mathbb {W}).\label{Def11}
\end{equation}
If $h\mapsto D_h F$ is bounded, then there exists a unique $DF\in
L^2(\mathbb{W},\mu_\mathbb{W};\mathbb{H})$ such that
\[
\langle DF,h\rangle _\mathbb{H}=D_hF\qquad \forall h\in \mathbb{H}.
\]
In this case, we shall write $F\in\mathscr{D}(D)$ and call $DF$ the
Malliavin gradient of~$F$ (cf. \cite{Nu}).

For $\ell\in\mathbb{S}_0$ and $x\in\mathbb{R}^d$, let $X^\ell
_t(x)=X^\ell_t$ solve the following SDE:
%
\begin{equation}
X^\ell_t=x+\int^t_0 b
\bigl(X^\ell_s\bigr)\,\mathrm{d}s+AW_{\ell_t}.\label{Ep1}
\end{equation}
Let $J^\ell_t:=J^\ell_t(x):=\nabla X^\ell_t(x)$ be the derivative matrix
of $X^\ell_t(x)$ with respect to the initial value $x$. It is easy to
see that
%
\begin{equation}
J^\ell_t=I+\int^t_0
\nabla b\bigl(X^\ell_s\bigr)\cdot J^\ell_s
\,\mathrm{d}s.\label{Ep2}
\end{equation}
Let $K^\ell_t$ be the inverse matrix of $J^\ell_t$. Then $K^\ell_t$ satisfies
%
\begin{equation}
K^\ell_t=I-\int^t_0K^\ell_s
\cdot\nabla b\bigl(X^\ell_s\bigr)\,\mathrm{d}s.\label{Ep3}
\end{equation}
Moreover, by definition (\ref{Def11}) and equation (\ref{Ep1}), it is
easy to see that $X^\ell_t(x)\in\mathscr{D}(D)$ and for any $h\in
\mathbb{H}$,
\[
D_hX^\ell_t=\int^t_0
\nabla b\bigl(X^\ell_s\bigr)D_hX^\ell_s
\,\mathrm{d}s+ Ah_{\ell_t}.
\]
The Malliavin covariance matrix is defined by
%
\begin{equation}
\bigl(\Sigma^\ell_t\bigr)_{ij}:=\bigl\langle D
\bigl(X^\ell_t\bigr)^i,D\bigl(X^\ell_t
\bigr)^j\bigr\rangle _\mathbb {H}.\label{ETR1}
\end{equation}

The following lemma provides an explicit expression of $\Sigma^\ell
_t$ in terms of $J^\ell_t$ (cf. \cite{Ku}),
which is crucial in the Malliavin's proof of H\"ormander's
hypoellipticity theorem.
%
%
\begin{lemma}
We have
%
\begin{equation}
\Sigma^\ell_t=J^\ell_t \biggl(\int
^t_0K^\ell_s AA^*
\bigl(K^\ell _s\bigr)^*\,\mathrm{d}\ell_s
\biggr) \bigl(J^\ell_t\bigr)^*,\label{Ep4}
\end{equation}
where $*$ denotes the transpose of a matrix.
\end{lemma}
\begin{pf}
For $\varepsilon\in(0,1)$, we define
%
\begin{equation}
\ell^\varepsilon_t:=\frac{1}{\varepsilon}\int^{t+\varepsilon
}_t
\ell_s\,\mathrm{d}s=\int^1_0
\ell_{\varepsilon s+t}\,\mathrm{d}s.\label{Def}
\end{equation}
Since $t\mapsto\ell_t$ is strictly increasing and right continuous,
it follows that for each $t\geq0$,
%
\begin{equation}
\ell^\varepsilon_t\downarrow\ell_t\qquad\mbox{as }\varepsilon\downarrow0.\label{ER1}
\end{equation}
Moreover, $t\mapsto\ell^\varepsilon_t$ is absolutely continuous and
strictly increasing.
Let $\gamma^\varepsilon$ be the inverse function of $\ell
^\varepsilon$, that is,
\[
\ell^\varepsilon_{\gamma^\varepsilon_t}=t,\qquad t\geq\ell ^\varepsilon_0
\quad\mbox{and}\quad \gamma^\varepsilon_{\ell^\varepsilon_t}=t,\qquad t\geq0.
\]
By definition, $\gamma^\varepsilon_t$ is also absolutely continuous
on $[\ell^\varepsilon_0,\infty)$.
Let $X^{\ell^\varepsilon}_t$ solve the following SDE:
\[
X^{\ell^\varepsilon}_t=x+\int^t_0 b
\bigl(X^{\ell^\varepsilon}_s\bigr)\,\mathrm{d} s+A(W_{\ell^\varepsilon_t}-W_{\ell^\varepsilon_0}).
\]
Let us now define
\[
Y^{\ell^\varepsilon}_t(x):=X^{\ell^\varepsilon}_{\gamma
^\varepsilon_t}(x),\qquad t\geq
\ell^\varepsilon_0.
\]
By the change of variables, one sees that
\[
Y^{\ell^\varepsilon}_t=x+\int^t_{\ell^\varepsilon_0} b
\bigl(Y^{\ell
^\varepsilon}_s\bigr)\dot\gamma^\varepsilon_s
\,\mathrm{d}s+A(W_t-W_{\ell
^\varepsilon_0}).
\]
It is well known that [cf. \cite{Nu}, p. 127, (2.60)]
\[
\bigl\langle DY^{\ell^\varepsilon}_t,DY^{\ell^\varepsilon}_t
\bigr\rangle _\mathbb{H}= \nabla Y^{\ell^\varepsilon}_t \biggl(
\int^t_{\ell^\varepsilon
_0}\bigl(\nabla Y^{\ell^\varepsilon}_s
\bigr)^{-1} AA^*\bigl(\bigl(\nabla Y^{\ell
^\varepsilon}_s
\bigr)^{-1}\bigr)^*\,\mathrm{d}s \biggr) \bigl(\nabla Y^{\ell^\varepsilon}_t
\bigr)^*.
\]
By the change of variables again, we obtain
%
\begin{eqnarray}\label{Ep5}
\qquad \bigl\langle DX^{\ell^\varepsilon}_t,DX^{\ell^\varepsilon}_t
\bigr\rangle _\mathbb{H} &=&\nabla X^{\ell^\varepsilon}_t \biggl(
\int^{\ell^\varepsilon
_t}_{\ell^\varepsilon_0}\bigl(\nabla Y^{\ell^\varepsilon}_s
\bigr)^{-1} AA^*\bigl(\bigl(\nabla Y^{\ell^\varepsilon}_s
\bigr)^{-1}\bigr)^*\,\mathrm{d}s \biggr) \bigl(\nabla X^{\ell^\varepsilon}_t
\bigr)^* \nonumber
\\
&=&\nabla X^{\ell^\varepsilon}_t \biggl(\int^t_0
\bigl(\nabla X^{\ell
^\varepsilon}_s\bigr)^{-1} AA^*\bigl(\bigl(
\nabla X^{\ell^\varepsilon}_s\bigr)^{-1}\bigr)^*\,\mathrm{d}\ell
^\varepsilon _s \biggr) \bigl(\nabla X^{\ell^\varepsilon}_t
\bigr)^*
\\
&=&J^{\ell^\varepsilon}_t \biggl(\int^t_0K^{\ell^\varepsilon
}_sAA^*
\bigl(K^{\ell^\varepsilon}_s\bigr)^*\,\mathrm{d}\ell^\varepsilon
_s \biggr) \bigl(J^{\ell^\varepsilon}_t\bigr)^*.\nonumber
\end{eqnarray}
From equation (\ref{Ep1}), it is easy to see that for each $t\geq0$
and $w\in\mathbb{W}$,
\[
\lim_{\varepsilon\downarrow0}\bigl|X^{\ell^\varepsilon}_t(w)-X^\ell
_t(w)\bigr|\leq C\lim_{\varepsilon\downarrow0}\bigl|W_{\ell^\varepsilon
_t}(w)-W_{\ell_t}(w)\bigr|=0.
\]
Thus, by equations (\ref{Ep2}) and (\ref{Ep3}), we also have
\[
\lim_{\varepsilon\downarrow0}\sup_{s\in[0,t]}\bigl|J^{\ell
^\varepsilon}_s(w)-J^\ell_s(w)\bigr|=0
\]
and
\[
\lim_{\varepsilon\downarrow0}\sup_{s\in[0,t]}\bigl|K^{\ell
^\varepsilon}_s(w)-K^\ell_s(w)\bigr|=0.
\]
Taking limits for both sides of (\ref{Ep5}), we obtain (\ref{Ep4})
(see \cite{Zh1}).
\end{pf}

The following lemma is a direct application of It\^o's formula (cf.
\cite{Pr}, p.~81, Theorem 33).
%
%
\begin{lemma}\label{Le2}
Let $V\dvtx \mathbb{R}^d\to\mathbb{M}^d$ be a $d\times d$-matrix valued
smooth function. We have
\begin{eqnarray*}
K^\ell_tV\bigl(X^\ell_t
\bigr)&=&V(x)+\int^t_0K^\ell_s(b
\cdot\nabla V-\nabla b\cdot V) \bigl(X^\ell_s\bigr)
\,\mathrm{d}s
\\
&&{} +\sum_{0<s\leq t}K^\ell_s
\bigl(V\bigl(X^\ell_s\bigr)-V\bigl(X^\ell
_{s-}\bigr)-\nabla V\bigl(X^\ell_{s-}\bigr)\cdot
\Delta X^\ell_s \bigr)
\\
&&{}+\int^t_0K^\ell_s
\cdot(\nabla V) \bigl(X^\ell_{s-}\bigr)\cdot A\,\mathrm{d}
W_{\ell_s},
\end{eqnarray*}
where $\Delta X^\ell_s:=X^\ell_s-X^\ell_{s-}=A(W_{\ell_s}-W_{\ell_{s-}})$.
\end{lemma}

We are now in a position to give the following.

\begin{pf*}{Proof of Theorem~\ref{Th1}}
By Lemma~\ref{EU1} and (\ref{Ep6}),
it is enough to prove that for each $\ell\in\mathbb{S}_0$, the law of
$X^\ell_t$ under $\mu_\mathbb{W}$ is absolutely continuous with
respect to the Lebesgue measure. By \cite{Nu}, page~97, Theorem 2.1.2,
it suffices to prove that $\Sigma^\ell_t$ is invertible. Since
$J^\ell_t$ is invertible, by (\ref{Ep4})
we only need to show that for any row vector $a\neq0\in\mathbb{R}^d$,
%
\begin{equation}
\int^t_0\bigl|aK^\ell_sA\bigr|^2
\,\mathrm{d}\ell_s>0.\label{Ep8}
\end{equation}
Suppose that
\[
\int^t_0\bigl|aK^\ell_sA\bigr|^2
\,\mathrm{d}\ell_s=\sum_{s\in(0,t]}\bigl|aK^\ell
_sA\bigr|^2\Delta\ell_s=0,
\]
then by Lemma~\ref{EU1} and the continuity of $s\mapsto|aK^\ell
_sA|$, we have
\[
aK^\ell_sA=0\qquad \forall s\in[0,t].
\]
Thus, by (\ref{Ep3}) we get
\[
0=aK^\ell_{t'}A=aA-\int^{t'}_0aK^\ell_s(
\nabla b) \bigl(X^\ell_s\bigr)A\,\mathrm{d} s\qquad \forall
t'\in[0,t],
\]
which in turn implies that
%
\begin{equation}\label{EE3}
aA=0
\end{equation}
and by the right continuity of $s\mapsto X^\ell_s$,
%
\begin{equation}
aK^\ell_sB_1\bigl(X^\ell_s
\bigr)A=aK^\ell_s(\nabla b) \bigl(X^\ell_s
\bigr)A=0\qquad \forall s\in[0,t].\label{EE2}
\end{equation}

Now we use the induction to prove that for each $n\in\mathbb{N}$,
%
\begin{equation}
aK^\ell_sB_n\bigl(X^\ell_s
\bigr)A=0\qquad \forall s\in[0,t].\label{EE22}
\end{equation}
Suppose that (\ref{EE22}) is true for some $n$. By Lemma~\ref{Le2},
we have
\[
K^\ell_tB_n\bigl(X^\ell_t
\bigr)=B_n(x)+\int^t_0K^\ell_sB_{n+1}
\bigl(X^\ell _s\bigr)\,\mathrm{d} s+M_t+V_t,
\]
where
\[
M_t:=\int^t_0K^\ell_s
\cdot(\nabla B_n) \bigl(X^\ell_{s-}\bigr)\cdot
A\,\mathrm{d} W_{\ell_s}
\]
and
\[
V_t:=\sum_{0<s\leq t}K^\ell_s
\bigl(B_n\bigl(X^\ell_s\bigr)-B_n
\bigl(X^\ell _{s-}\bigr)-(\nabla B_n)
\bigl(X^\ell_{s-}\bigr)\cdot\Delta(AW_{\ell_s})
\bigr).
\]
Thus, by (\ref{EE22}) we have
%
\begin{equation}
\int^{t'}_0aK^\ell_sB_{n+1}
\bigl(X^\ell_s\bigr)A\,\mathrm{d}s+aM_{t'}A+aV_{t'}A=0\qquad\forall t'\in[0,t].\label{EW1}
\end{equation}
By the inductive assumption (\ref{EE22}), we have
\[
aK^\ell_sB_n\bigl(X^\ell_s
\bigr)A=aK^\ell_sB_n\bigl(X^\ell_{s-}
\bigr)A=0.
\]
Hence,
\begin{eqnarray*}
aV_{t'}A&=&-\sum_{0<s\leq t'}aK^\ell_s
\cdot(\nabla B_n) \bigl(X^\ell _{s-}\bigr)\cdot
\Delta(AW_{\ell_s})\cdot A
\\
&=&-\int^{t'}_0aK^\ell_s
\cdot(\nabla B_n) \bigl(X^\ell_{s-}\bigr)\cdot
A\,\mathrm{d} W_{\ell_s}\cdot A=-aM_{t'}A,
\end{eqnarray*}
which together with (\ref{EW1}) implies that
\[
aK^\ell_sB_{n+1}\bigl(X^\ell_s
\bigr)A=0\qquad \forall s\in[0,t].
\]
The assertion (\ref{EE22}) is thus proved. Combining (\ref{EE3}) and
(\ref{EE22}) and by letting $s\to0$, we obtain
\[
aA=aB_1(x)A=\cdots=aB_n(x)A=0,
\]
which is contrary to (\ref{eqHHN}). The proof is thus complete.
\end{pf*}

\section{Proof of Theorem \texorpdfstring{\protect\ref{Th2}}{1.3}}\label{sec3}
\subsection{Norris' type lemma}\label{sec3.1}
In this section, we use the following filtration:
\[
\mathscr{F}_t:=\sigma\{W_{S_s},S_s\dvtx  s\leq t
\}.
\]
Clearly, for $t>s$, $W_{S_t}-W_{S_s}$ and $S_t-S_s$ are independent of
$\mathscr{F}_s$.

Let us first prove the following estimate of exponential type about the
subordinator $S_t$.
%
%
\begin{lemma}\label{Le1}
Let $f_t\dvtx \mathbb{R}_+\to\mathbb{R}_+$ be a bounded continuous
nonnegative \mbox{$\mathscr{F}_t$-}adapted process. For any $\varepsilon,\delta>0$, we have
\[
P \biggl\{\int^t_0f_s
\,\mathrm{d}S_s\leq\varepsilon;\int^t_0f_s
\,\mathrm{d} s>\delta \biggr\}\leq\mathrm{e}^{1-\phi(1/\varepsilon)\delta},
\]
where
\[
\phi(\lambda):=\frac{\lambda}{2}\int^{(\log2)/(\lambda\|f\|
_\infty)}_0 u
\nu_S(\mathrm{d}u),\qquad \lambda>0,
\]
and $\nu_S$ is the L\'evy measure of the subordinator $S_t$.
\end{lemma}
\begin{pf}
For $\lambda>0$, set
\[
g^\lambda_s:=\int^\infty_0
\bigl(1-\mathrm{e}^{-\lambda f_su}\bigr)\nu _S(\mathrm{d}u)
\]
and
\[
M^\lambda_t:=-\lambda\int^t_0f_s
\,\mathrm{d}S_s+\int^t_0g^\lambda
_s\,\mathrm{d}s.
\]
Let $\mu(t,\mathrm{d}u)$ be the Poisson random measure associated with
$S_t$, that is,
\[
\mu(t,U):=\sum_{s\leq t}1_U(\Delta
S_s),\qquad U\in\mathscr{B}(\mathbb{R}_+).
\]
Let $\tilde\mu(t,\mathrm{d}u)$ be the compensated Poisson random measure
of $\mu(t,\mathrm{d}u)$, that is,
\[
\tilde\mu(t,\mathrm{d}u)=\mu(t,\mathrm{d}u)-t\nu_S(\mathrm{d}u).
\]
Then we can write
\[
\int^t_0f_s
\,\mathrm{d}S_s=\int^t_0\!\int^\infty_0f_su\mu (\mathrm{d}s,
\mathrm{d}u).
\]
By It\^o's formula, we have
\[
\mathrm{e}^{M^\lambda_t}=1+\int^t_0\!\int
^\infty_0\mathrm {e}^{M^\lambda_{s-}}\bigl[
\mathrm{e}^{-\lambda f_su}-1\bigr]\tilde\mu(\mathrm{d} s,\mathrm{d}u).
\]
Since for $x>0$,
\[
1-\mathrm{e}^{-x}\leq1\wedge x,
\]
we have
\[
g^\lambda_s\leq\int^\infty_0
\bigl(1\wedge\bigl(\lambda\|f\|_\infty u\bigr)\bigr)\nu _S(
\mathrm{d}u)
\]
and
\[
M^\lambda_t\leq\int^t_0g^\lambda_s
\,\mathrm{d}s\leq t\int^\infty _0\bigl(1\wedge\bigl(
\lambda\|f\|_\infty u\bigr)\bigr)\nu_S(\mathrm{d}u).
\]
Hence, for any $\lambda>0$ and $t>0$,
\[
\mathbb{E}\mathrm{e}^{M^\lambda_t}=1.
\]
On the other hand, since for any $\kappa\in(0,1)$ and $0\leq x\leq
-\log k$,
\[
1-\mathrm{e}^{-x}\geq\kappa x,
\]
we have
\begin{eqnarray*}
g^\lambda_s&\geq&\int^{(\log2)/(\lambda\|f\|_\infty
)}_0\bigl(1-\mathrm{e}^{-\lambda f_su}\bigr)\nu_S(\mathrm{d}u)
\\
& \geq&
\frac{\lambda f_s}{2}\int^{(\log2)/(\lambda\|f\|_\infty
)}_0 u \nu_S(\mathrm{d}u) =\phi(\lambda) f_s.
\end{eqnarray*}
Thus,
\begin{eqnarray*}
\biggl\{\int^t_0f_s
\,\mathrm{d}S_s\leq\varepsilon;\int^t_0f_s
\,\mathrm{d} s>\delta \biggr\} &\subset& \biggl\{\mathrm{e}^{M^\lambda_t}\geq
\mathrm{e}^{-\lambda
\varepsilon+\int^t_0g^\lambda_s\,\mathrm{d}s};\int^t_0g^\lambda
_s\,\mathrm{d} s>\phi(\lambda)\delta \biggr\}
\\
&\subset& \bigl\{\mathrm{e}^{M^\lambda_t}\geq\mathrm{e}^{-\lambda
\varepsilon+\phi(\lambda)\delta} \bigr\},
\end{eqnarray*}
which then implies the result by Chebyshev's inequality and letting
$\lambda=\frac{1}{\varepsilon}$.
\end{pf}

Let $N(t,\mathrm{d}y)$ be the Poisson random measure associated with
$L_t=W_{S_t}$, that is,
\[
N(t,\Gamma)=\sum_{s\leq t}1_\Gamma(L_s-L_{s-}),\qquad
\Gamma\in \mathscr{B}\bigl(\mathbb{R}^d\bigr).
\]
Let $\widetilde N(t,\mathrm{d}y)$ be the compensated Poisson random
measure of
$N(t,\mathrm{d}y)$, that is,
\[
\widetilde N(t,\mathrm{d}y)=N(t,\mathrm{d}y)-t\nu_L(\mathrm{d}y),
\]
where $\nu_L$ is the L\'evy measure of $L_t$ given by (\ref{EW2}). By
L\'evy--It\^o's decomposition (cf. \cite{Ap}), we have
%
\begin{equation}
L_t=W_{S_t}=\int_{|y|\leq1}y\widetilde N(t,
\mathrm{d}y)+\int_{|y|>1}yN(t,\mathrm{d} y).\label{EW6}
\end{equation}

We recall the following result about the exponential estimate of
discontinuous martingales (cf. \cite{Ca}, Lemma 1).
%
%
\begin{lemma}\label{Lemma1}
Let $f_t(y)$ be a bounded $\mathscr{F}_t$-predictable process with
bound $A$.
Then for any $\delta,\rho>0$, we have
\begin{eqnarray*}
&&P \biggl\{\sup_{t\in[0,T]}\biggl\llvert \int^t_0\!\int_{\mathbb
{R}^d}f_s(y)\widetilde N(\mathrm{d}s,
\mathrm{d}y)\biggr\rrvert \geq\delta, \int^T_0\!\!\int_{\mathbb{R}^d}\bigl|f_s(y)\bigr|^2
\nu_L(\mathrm {d}y)\,\mathrm{d}s<\rho \biggr\}
\\
&&\qquad \leq2\exp \biggl(-\frac{\delta^2}{2(A\delta
+\rho)} \biggr).
\end{eqnarray*}
\end{lemma}
The following lemma is contained in the proof of Norris' lemma (cf.
\cite{Nu}, p.~137).
%
%
\begin{lemma}\label{Le4}
For $T>0$, let $f$ be a bounded measurable $\mathbb{R}^d$-valued
function on $[0,T]$. Assume that for some $\varepsilon<T$ and $x\in
\mathbb{R}^d$,
%
\begin{equation}
\int^T_0\biggl\llvert x+\int
^t_0f_s\,\mathrm{d}s\biggr\rrvert
^2\,\mathrm{d}t\leq \varepsilon ^3.\label{ET4}
\end{equation}
Then we have
\[
\sup_{t\in[0,T]}\biggl\llvert \int^t_0f_s
\,\mathrm{d}s\biggr\rrvert \leq2\bigl(1+\|f\| _\infty\bigr)\varepsilon.
\]
\end{lemma}
\begin{pf}
By (\ref{ET4}) and Chebyshev's inequality, we have
\[
\mathrm{Leb} \biggl\{t\in[0,T]\dvtx  \biggl\llvert x+\int^t_0f_s
\,\mathrm{d}s\biggr\rrvert \geq\varepsilon \biggr\}\leq\varepsilon<T.
\]
Thus, for each $t\in[0,T]$, there exits an $s\in[0,T]$ such that
\[
|s-t|\leq\varepsilon\quad\mbox{and}\quad \biggl\llvert x+\int^s_0f_r
\,\mathrm{d} r\biggr\rrvert <\varepsilon.
\]
Consequently, for such $t,s$,
\[
\biggl\llvert x+\int^t_0f_r
\,\mathrm{d}r\biggr\rrvert \leq\biggl\llvert x+\int^s_0f_r
\,\mathrm{d} r\biggr\rrvert +\biggl\llvert \int^t_sf_r
\,\mathrm{d}r\biggr\rrvert \leq\varepsilon+\varepsilon\|f\|_\infty.
\]
In particular,
\[
|x|\leq\varepsilon+\varepsilon\|f\|_\infty,
\]
hence,
\[
\biggl\llvert \int^t_0f_s
\,\mathrm{d}s\biggr\rrvert \leq|x|+\biggl\llvert x+\int^t_0f_s
\,\mathrm{d} s\biggr\rrvert \leq2\bigl(\varepsilon+\varepsilon\|f\|_\infty\bigr).
\]
The proof is finished.
\end{pf}

We now prove the following Norris' type lemma (cf. \cite{Nu,Ca}).
%
%
\begin{lemma}\label{Lemma2}
Let $Y_t=y+\int^t_0\beta_s\,\mathrm{d}s$ be an $\mathbb{R}^d$-valued
process, where $\beta_t$ takes the following form:
\[
\beta_t=\beta_0+\int^t_0
\gamma_s\,\mathrm{d}s+\int^t_0\!\int_{\mathbb{R}^d}g_s(y)\widetilde N(\mathrm{d}s,\mathrm{d}y),
\]
where $\gamma_t$ and $g_t(y)$ are two $\mathscr{F}_t$-predictable
$\mathbb{R}^d$-valued processes.
Suppose that for some nonrandom constants $C_1,C_2\geq1$ and all
$s\geq0, y\in\mathbb{R}^d$,
%
\begin{equation}
|\beta_s|+|\gamma_s|\leq C_1,\qquad \bigl|g_s(y)\bigr|\leq C_2\bigl(1\wedge|y|\bigr).
\end{equation}
Then for any $\delta\in(0,\frac{1}{3})$, there exists $\varepsilon
_0=\varepsilon_0(C_1, C_2,\nu_L,\delta)\in(0,1)$
such that for all $T\in(0,1)$ and $\varepsilon\in(0,T^3\wedge
\varepsilon_0)$,
%
\begin{equation}
P \biggl\{\int^T_0|Y_s|^2
\,\mathrm{d}s<\varepsilon, \int^T_0|\beta
_s|^2\,\mathrm{d}s\geq9C^2_1
\varepsilon^{\delta} \biggr\}\leq 2\exp \biggl\{-\frac{\varepsilon^{\delta-(1/3)}}{9C_1} \biggr
\}.\label{EW5}
\end{equation}
\end{lemma}
\begin{pf}
Let us define
\[
h_t:=\int^t_0
\beta_s\,\mathrm{d}s,\qquad M_t:=\int^t_0\!\int_{\mathbb{R}^d}\bigl\langle h_s, g_s(y)
\bigr\rangle _{\mathbb{R}^d}\widetilde N(\mathrm{d}s,\mathrm{d}y)
\]
and
\begin{eqnarray*}
E_1&:=& \biggl\{\int^T_0|Y_s|^2
\,\mathrm{d}s<\varepsilon \biggr\},\qquad E_2:= \Bigl\{\sup
_{t\in[0,T]}|h_t|\leq2(1+C_1)
\varepsilon^{1/3} \Bigr\},
\\
E_3&:=& \bigl\{\langle M\rangle _T \leq
C_3\varepsilon^{2/3} \bigr\},\qquad E_4:= \Bigl\{\sup_{t\in[0,T]}|M_t|\leq
\varepsilon ^\delta \Bigr\},
\\
E_5&:=& \biggl\{\int^T_0|
\beta_s|^2\,\mathrm{d}s<9C_1^2
\varepsilon ^\delta \biggr\},
\end{eqnarray*}
where $C_3$ is determined below.

First of all, by Lemma~\ref{Le4}, one sees that for $\varepsilon<T^3$,
%
\begin{equation}
E_1\subset E_2\subset E_3,\label{EK1}
\end{equation}
where the second inclusion is due to
\begin{eqnarray*}
\langle M\rangle _T&=&\int^T_0\!\!\int_{\mathbb{R}^d}\bigl|\bigl\langle h_s, g_s(y)
\bigr\rangle _{\mathbb
{R}^d}\bigr|^2\nu_L(\mathrm{d}y)
\,\mathrm{d}s
\\
&\leq&4(1+C_1)^2C_2^2 \biggl(\int
_{\mathbb{R}^d}1\wedge|y|^2\nu _L(
\mathrm{d}y) \biggr)\varepsilon^{2/3}=:C_3\varepsilon
^{2/3}.
\end{eqnarray*}
On the other hand, by the integration by parts formula, we have
\[
\int^T_0|\beta_t|^2
\,\mathrm{d}t=\int^T_0\langle
\beta_t, \mathrm {d}h_t\rangle _{\mathbb{R}^d} =\langle
\beta_T, h_T\rangle _{\mathbb{R}^d}-\int
^T_0\langle h_t,
\gamma_t\rangle _{\mathbb
{R}^d}\,\mathrm{d}t-M_T.
\]
From this, one sees that on $E_2\cap E_4$,
\begin{eqnarray*}
\int^T_0|\beta_t|^2
\,\mathrm{d}t&\leq&2C_1(1+C_1)\varepsilon^{1/3}(1+T)+
\varepsilon^\delta
\\
&\leq&\bigl(4C_1(1+C_1)+1\bigr)\varepsilon^\delta
\leq9C^2_1\varepsilon ^\delta.
\end{eqnarray*}
This means that
\[
E_2\cap E_4\subset E_5,
\]
which together with (\ref{EK1}) gives
\[
E_1\cap E_5^c\subset E_1\cap
E^c_4\subset E_2\cap E_3\cap
E^c_4.
\]
Thus, by Lemma~\ref{Lemma1} we have
\[
P\bigl(E_1\cap E_5^c\bigr)\leq2\exp
\biggl(-\frac{\varepsilon^{2\delta
}}{2(2(1+C_1)\varepsilon^{(1/3)+\delta}+C_3\varepsilon
^{2/3})} \biggr)
\]
and (\ref{EW5}) follows by choosing $\varepsilon_0$ with
$C_3\varepsilon^{(1/3)-\delta}_0=1$.
\end{pf}

Below we set
%
\begin{equation}\label{ETR2}
\Sigma_t:=\Sigma^\ell_t|_{\ell=S},\qquad
K_t:=K^\ell_t|_{\ell=S},\qquad
J_t:=J^\ell_t|_{\ell=S}.
\end{equation}
The following lemma is a key step for proving the smoothness of $p_t(x,y)$.
%
%
\begin{lemma}\label{Le35}
Let $\theta\in(0,\frac{1}{2})$ be given in (\ref{Con2}). Under
(\ref{eqUH1}) and (\ref{Con2}), for any $p>1$, there exist
$C_0=C_0(p,\theta)>0$
and $C_1=C_1(p,\theta)>0$ such that for all $t\in(0,1)$ and
$\varepsilon\in(0,C_0t^{8/\theta})$,
%
\begin{equation}
\sup_{|a|=1}P \biggl\{\int^t_0|aK_sA|^2
\,\mathrm{d}S_s\leq\varepsilon \biggr\}\leq C_1
\varepsilon^p.\label{EW4}
\end{equation}
\end{lemma}
\begin{pf}
By Lemma~\ref{Le1} and (\ref{Con2}), for the given $\theta$ in (\ref{Con2}),
there exists an $\varepsilon_0=\varepsilon_0(\theta)>0$ such that
for all $\varepsilon\in(0,\varepsilon_0)$ and $t\in(0,1)$,
%
\begin{eqnarray}\label{EW3}
&&P \biggl\{\int^t_0|aK_sA|^2
\,\mathrm{d}S_s\leq\varepsilon \biggr\}
\nonumber
\\
&&\qquad\leq P \biggl\{\int^t_0|aK_sA|^2
\,\mathrm{d}S_s\leq\varepsilon, \int^t_0|aK_sA|^2
\,\mathrm{d}s\geq\varepsilon^\theta \biggr\}
\nonumber
\\
&&\qquad\quad{}+P \biggl\{\int^t_0|aK_sA|^2
\,\mathrm{d}s<\varepsilon ^\theta \biggr\}
\\
&&\qquad\leq\exp \biggl\{1-\frac{1}{2\varepsilon^{1-\theta}}\int^{C\varepsilon}_0u
\nu_S(\mathrm{d}u) \biggr\} +P \biggl\{\int^t_0|aK_sA|^2
\,\mathrm{d}s<\varepsilon^\theta \biggr\}
\nonumber
\\
&&\qquad\leq\exp \bigl\{1-\varepsilon^{-\theta/2} \bigr\}+P \biggl\{\int
^t_0|aK_sA|^2
\,\mathrm{d}s<\varepsilon^\theta \biggr\}.\nonumber
\end{eqnarray}
Notice that by (\ref{EW6}),
\[
X_t=x+\int^t_0b(X_s)
\,\mathrm{d}s+\int_{|y|\leq1}Ay\widetilde N(t,\mathrm {d}y)+\int
_{|y|>1}Ay N(t,\mathrm{d}y).
\]
If we set $Y_t:=aK_tA$ and
\begin{eqnarray*}
\beta_t&:=&aK_t\nabla b(X_t)A,\qquad
g_t(y):=aK_t\bigl(\nabla b(X_{t-}+Ay)-\nabla
b(X_{t-})\bigr)A,
\\
\gamma_t&:=&\int_{\mathbb{R}^d}aK_t \bigl(
\nabla b(X_t+Ay)-\nabla b(X_t)-1_{|y|\leq1}Ay\cdot
\nabla^2 b(X_t) \bigr)A\nu_L(\mathrm {d}y)
\\
&&{}+aK_tB_2(X_t)A,
\end{eqnarray*}
then by equation (\ref{Ep3}) and It\^o's formula, one sees that
$Y_t=aA+\int^t_0\beta_s\,\mathrm{d}s$ and
\[
\beta_t=a\nabla b(x)A+\int^t_0
\gamma_s\,\mathrm{d}s+\int^t_0\!\int_{\mathbb{R}^d}g_s(y)\widetilde N(\mathrm{d}s,\mathrm{d}y).
\]
By the assumptions, it is easy to see that
\[
|\beta_t|+|\gamma_t|\leq C_1\bigl(\|\nabla b
\|_\infty,\bigl\|\nabla^3 b\bigr\| _\infty, \|A\|\bigr)
\]
and
\[
\bigl|g_t(y)\bigr|\leq C_2\bigl(\|\nabla b\|_\infty,\bigl\|
\nabla^2 b\bigr\|_\infty, \|A\| \bigr) \bigl(1\wedge|y|\bigr).
\]
Fix $\delta\in(0,\frac{1}{3})$. Define now
\[
E^\varepsilon_t:= \biggl\{\int^t_0|aK_sA|^2
\,\mathrm{d}s<\varepsilon ^\theta \biggr\},\qquad F^\varepsilon_t:=
\biggl\{\int^t_0\bigl|aK_s\nabla
b(X_s)A\bigr|^2\,\mathrm{d} s<9C_1^2
\varepsilon^{\theta\delta} \biggr\}.
\]
Then, by Lemma~\ref{Lemma2}, there is an $\varepsilon_0\in(0,1)$
such that for all $t\in(0,1)$ and
$\varepsilon\in(0,t^3\wedge\varepsilon_0)$,
\begin{eqnarray*}
P\bigl(E^\varepsilon_t\bigr)&=&P\bigl(E^\varepsilon_t
\cap\bigl(F^\varepsilon _t\bigr)^c\bigr)+P
\bigl(E^\varepsilon_t\cap F^\varepsilon_t\bigr)
\\
&\leq&2\exp\bigl\{-\varepsilon^{\theta(\delta-(1/3))}/(9C_1)\bigr\} +P
\bigl(E^\varepsilon_t\cap F^\varepsilon_t\bigr).
\end{eqnarray*}
Define
\[
\tau:=\inf \bigl\{s\geq0\dvtx  |K_s-I|\geq\tfrac{1}{2} \bigr\}
\wedge t.
\]
Then
\[
P\bigl(E^\varepsilon_t\cap F^\varepsilon_t\bigr)
\leq P\bigl(E^\varepsilon_t\cap F^\varepsilon_t
\cap\bigl\{\tau\geq\varepsilon^{\delta\theta/2}\bigr\} \bigr)+P\bigl(\tau<
\varepsilon^{\delta\theta/2}\bigr).
\]
By Chebyshev's inequality, we have for any $p>1$,
\begin{eqnarray*}
P\bigl(\tau<\varepsilon^{\delta\theta/2}\bigr)&\leq& P \biggl\{\sup
_{s\in
(0,\varepsilon^{\delta\theta/2}\wedge t)}|K_s-I|\geq\frac
{1}{2} \biggr\}
\\
&\leq& 2^p\mathbb{E} \Bigl(\sup_{s\in(0,\varepsilon^{\delta\theta
/2}\wedge t)}|K_s-I|^p
\Bigr)
\\
&\leq& C\bigl(\varepsilon^{\delta\theta/2}\wedge t\bigr)^p
\end{eqnarray*}
and by (\ref{eqUH1}),
\begin{eqnarray*}
E^\varepsilon_t\cap F^\varepsilon_t&\subset&
\biggl\{\int^t_0\bigl(|aK_sA|^2+\bigl|aK_s
\nabla b(X_s)A\bigr|^2\bigr)\,\mathrm{d}s<\varepsilon^\theta
+9C^2_1\varepsilon^{\delta\theta} \biggr\}
\\
&\subset& \biggl\{\int^t_0
\frac{|aK_sA|^2+|aK_s\nabla
b(X_s)A|^2}{|aK_s|^2}|aK_s|^2\,\mathrm{d}s<
\bigl(1+9C^2_1\bigr)\varepsilon ^{\delta
\theta} \biggr\}
\\
&\subset& \biggl\{c_1\int^t_0|aK_s|^2
\,\mathrm{d}s<\bigl(1+9C^2_1\bigr)\varepsilon
^{\delta\theta} \biggr\}.
\end{eqnarray*}
Since on $\{\tau\geq\varepsilon^{\delta\theta/2}\}$,
\[
|aK_s|\geq1-|K_s-I|\geq\tfrac{1}{2},\qquad
|a|=1,\qquad
s\in \bigl[0,\varepsilon^{\delta\theta/2}\wedge t\bigr],
\]
it is easy to see that for any $\varepsilon<t^{2/\delta
\theta}\wedge(\frac{c_1}{4(1+9C^2_1)})^{2/(\delta\theta)}$,
\[
E^\varepsilon_t\cap F^\varepsilon_t\cap\bigl\{
\tau\geq\varepsilon ^{\delta\theta/2}\bigr\}\subset \bigl\{c_1\bigl(
\varepsilon^{\delta\theta
/2}\wedge t\bigr)/4<\bigl(1+9C^2_1
\bigr)\varepsilon^{\delta\theta} \bigr\} =\varnothing.
\]
Hence, for any $p>1$,
if one takes $\delta=\frac{1}{4}$ and $C_0=C_0(\varepsilon
_0,p,\theta,c_1)$ being small enough, then
for all $t\in(0,1)$ and $\varepsilon\in(0,C_0t^{8/\theta})$,
\[
P\bigl(E^\varepsilon_t\bigr)\leq C\varepsilon^{\theta p/8},
\]
which together with (\ref{EW3}) yields (\ref{EW4}) by resetting
$p=\frac{8p'}{\theta}$.
\end{pf}

\subsection{$S_t$ has finite moments of all orders}\label{sec3.2}

In this subsection, we suppose that $S_t$ has finite moments of all
orders and
$b\in C^\infty(\mathbb{R}^d)$ has bounded derivatives of all orders.
The following lemma is standard.
%
%
\begin{lemma}\label{Le3}
For any $m,k\in\{0\}\cup\mathbb{N}$ with $m+k\geq1$ and $p\geq1$,
we have
%
\begin{equation}
\sup_{x\in\mathbb{R}^d}\sup_{t\in[0,1]}\mathbb{E} \bigl(\bigl\|
D^m\nabla^k X^\ell_t(x)
\bigr\|_{\mathbb{H}^{\otimes^m}}^p |_{\ell=S} \bigr)<+\infty.\label{ET6}
\end{equation}
\end{lemma}
\begin{pf}Noticing that
\[
DX_t^\ell(x)=\int^t_0
\nabla b\bigl(X^\ell_s(x)\bigr)DX_s^\ell(x)
\,\mathrm{d}s+\cdot\wedge\ell_t
\]
and
\[
\nabla X^\ell_t(x)=I+\int^t_0
\nabla b\bigl(X^\ell_s(x)\bigr)\nabla X^\ell
_s(x)\,\mathrm{d}s,
\]
we have
\[
\bigl\|DX_t^\ell(x)\bigr\|_\mathbb{H}\leq\|\nabla b
\|_\infty\int^t_0\bigl\|
DX_s^\ell(x)\bigr\|_\mathbb{H}\,\mathrm{d}s+
\ell_t^{1/2}
\]
and
\[
\bigl|\nabla X^\ell_t(x)\bigr|\leq1+\|\nabla b\|_\infty
\int^t_0\bigl|\nabla X^\ell
_s(x)\bigr|\,\mathrm{d}s.
\]
By Gronwall's inequality, we obtain
\[
\bigl\|DX_t^\ell(x)\bigr\|_\mathbb{H}\leq
\ell_t^{1/2}+\mathrm{e}^{\|\nabla
b\|_\infty t}\int
^t_0\ell_s^{1/2}
\,\mathrm{d}s
\]
and
\[
\bigl\|\nabla X_t^\ell(x)\bigr\|_\mathbb{H}\leq
\mathrm{e}^{\|\nabla b\|
_\infty t}.
\]
Hence, for any $p\geq1$,
\[
\mathbb{E} \bigl(\bigl\|DX^\ell_t\bigr\|_{\mathbb{H}}^p
|_{\ell=S} \bigr)\leq C\mathbb{E}|S_t|^{p/2}+C\int
^t_0\mathbb{E}|S_s|^{p/2}
\,\mathrm{d} s<+\infty.
\]
Thus, we obtain (\ref{ET6}) for $m+k=1$. For the general $m$ and $k$,
it follows by similar calculations and the induction.
\end{pf}

We recall the following main criterion in the Malliavin calculus that a
random vector admits a smooth density (cf. \cite{Nu}, pp.~100--103).
%
%
\begin{proposition}\label{Pro2}
Let $F=(F^1,\ldots, F^d)$ be a smooth Wiener functional and $(\Sigma
_F)_{ij}:=\langle DF^i, DF^j\rangle _\mathbb{H}$ be the Malliavin covariance
matrix. We assume that for all $p\geq2$,
\[
\mathbb{E}\bigl[(\det\Sigma_F)^{-p}\bigr]<\infty.
\]
Let $G$ be another smooth Wiener functional and $\varphi\in C^\infty
_b(\mathbb{R}^d)$. Then for any multi-index $\alpha=(\alpha_1,\ldots,\alpha_m)\in\{1,2,\ldots,d\}^m$,
\[
\mathbb{E}\bigl[\partial_\alpha\varphi(F)G\bigr]=\mathbb{E}\bigl[
\varphi (F)H_{\alpha}(F,G)\bigr],
\]
where $\partial_\alpha=\partial_{\alpha_1}\cdots\partial
_{\alpha_m}$, and $H_{\alpha}(F,G)$ are recursively defined by
\begin{eqnarray*}
H_{(i)}(F,G)&:=&\sum_{j}D^* \bigl(G
\bigl(\Sigma^{-1}_F\bigr)_{ij}DF^j
\bigr),
\\
H_\alpha(F,G)&:=&H_{(\alpha_m)}\bigl(F,H_{(\alpha_1,\ldots,\alpha_{m-1})}(F,G)\bigr).
\end{eqnarray*}
As a consequence, for any $p\geq1$, there exist $p_1,p_2,p_3>1$ and
$n_1,n_2\in\mathbb{N}$ such that
\[
\bigl\|H_\alpha(F,G)\bigr\|_p\leq C\bigl\|(\det\Sigma_F)^{-1}
\bigr\|^{n_1}_{p_1}\|DF\| _{m,p_2}^{n_2}\|G
\|_{m,p_3}.
\]
In particular, the law of $F$ possesses an infinitely differentiable
density $\rho\in\mathcal{S}(\mathbb{R}^d)$, the space of Schwartz
rapidly decreasing functions.
\end{proposition}

Now we can prove the following gradient estimate.
%
%
\begin{theorem}\label{Th4}
Under (\ref{eqUH1}) and (\ref{Con2}), for any $k,m\in\{0\}\cup
\mathbb{N}$ with \mbox{$k+m\geq1$}, there are $\gamma_{k,m}>0$
and $C=C(k,m)>0$ such that for all $f\in C^\infty_b(\mathbb{R}^d)$
and $t\in(0,1)$,
%
\begin{equation}
\sup_{x\in\mathbb{R}^d} \bigl|\nabla^k\mathbb{E} \bigl( \bigl(
\nabla ^mf\bigr) \bigl(X_t(x)\bigr) \bigr) \bigr|\leq C\|f
\|_\infty t^{-\gamma_{k,m}}.\label{ET7}
\end{equation}
\end{theorem}
\begin{pf} We first prove that there exists a constant $\gamma>0$
such that
for any $p\geq1$, some $C=C(p)>0$ and all $t\in(0,1)$,
%
\begin{equation}
\bigl\|(\det\Sigma_t)^{-1}\bigr\|_p\leq
Ct^{-\gamma},\label{Es1}
\end{equation}
which, by (\ref{Ep4}), is equivalent to prove that
\[
\biggl\llVert \det \biggl(\int^t_0K_s
AA^*K_s^*\,\mathrm{d}S_s \biggr)^{-1}\biggr
\rrVert _p\leq Ct^{-\gamma}.
\]
Since the determinant of a matrix is greater than \mbox{$d$-}times its
smallest eigenvalue, that is,
\[
\biggl(\inf_{|a|=1}\int^t_0|aK_s
A|^2\,\mathrm{d}S_s \biggr)^d\leq \det \biggl(
\int^t_0K_s AA^*K_s^*
\,\mathrm{d}S_s \biggr),
\]
it suffices to prove that for some $\gamma'>0$,
\[
\biggl\llVert \biggl(\inf_{|a|=1}\int^t_0|aK_s
A|^2\,\mathrm{d}S_s \biggr)^{-1}\biggr\rrVert
_p\leq Ct^{-\gamma'},
\]
which will follow by showing that for all $p\geq1$ and $\varepsilon
\in(0,C_pt^{\gamma'})$,
\[
P \biggl\{\inf_{|a|=1}\int^t_0|aK_s
A|^2\,\mathrm{d}S_s\leq\varepsilon \biggr\}\leq C
\varepsilon^p.
\]
Since $S_t$ has finite moments of all orders, this estimate follows by
(\ref{EW4}) and a compact argument
(see \cite{Nu}, p. 133, Lemma 2.3.1, for more details).

Next, by the chain rule, we have
\begin{eqnarray*}
&&\nabla^k\mathbb{E} \bigl(\bigl(\nabla^mf\bigr)
\bigl(X_t(x)\bigr) \bigr)
\\
&&\qquad  =\sum_{j=1}^k
\mathbb{E} \bigl(\bigl(\nabla^{m+j}f\bigr) \bigl(X_t(x)
\bigr)G_j\bigl(\nabla X_t(x),\ldots,\nabla^k
X_t(x)\bigr) \bigr)
\\
&&\qquad=\sum_{j=1}^k\mathbb{E} \bigl(
\mathbb{E} \bigl(\bigl(\nabla ^{m+j}f\bigr) \bigl(X^\ell_t(x)
\bigr)G_j\bigl(\nabla X^\ell_t(x),\ldots,
\nabla^k X^\ell _t(x)\bigr) \bigr)
|_{\ell=S} \bigr),
\end{eqnarray*}
where $\{G_j, j=1,\ldots,k\}$ are real polynomial functions.
By Proposition~\ref{Pro2}, Lemma~\ref{Le3} and H\"older's inequality,
there exist integer $n$ and $p>1$, $C>0$ such that for all $t\in(0,1)$,
\begin{eqnarray*}
\bigl|\nabla^k\mathbb{E} \bigl(\bigl(\nabla^mf\bigr)
\bigl(X_t(x)\bigr) \bigr) \bigr| &\leq& C\|f\|_\infty\mathbb{E}
\bigl(\bigl\|\bigl(\det\Sigma^\ell_t\bigr)^{-1}\bigr\|
^n_p |_{\ell=S} \bigr)
\\
&\leq& C\|f\|_\infty\bigl\|(\det\Sigma_t)^{-1}
\bigr\|^n_{np}.
\end{eqnarray*}
Estimate (\ref{ET7}) now follows by (\ref{Es1}).
\end{pf}

\subsection{Without the finiteness assumption of moments}\label{sec3.3}
Let $S'_t$ be a subordinator with L\'evy measure $1_{(0,1)}(u)\nu
_S(\mathrm{d}u)$ and independent of
$(W_t)_{t\geq0}$. Let $ p'_t(x,y)$ be the distributional density of $X'_t(x)$,
where $X'_t(x)$ solves the following SDE:
\[
X'_t(x)=x+\int^t_0b
\bigl(X'_s(x)\bigr)\,\mathrm{d}s+AW_{S'_t}.
\]
Let us write
\[
\mathcal{P}'_t f(x):=\mathbb{E}f\bigl(X'_t(x)
\bigr)=\int_{\mathbb{R}^d}f(y) p'_t(x,y)
\,\mathrm{d}y.
\]

We first prepare two simple lemmas for later use.
%
%
\begin{lemma}
Let $f\in C^\infty_b(\mathbb{R}^d)$. For any $m\in\mathbb{N}$,
there exists a constant \mbox{$C_{m,b}\geq1$}
such that for all $x\in\mathbb{R}^d$ and $t\in[0,1]$,
%
\begin{equation}
\bigl|\nabla^m\mathcal{P}'_tf(x)\bigr|\leq
C_{m,b}\sum_{k=1}^m\mathcal
{P}'_t\bigl|\nabla^kf\bigr|(x).\label{ET5}
\end{equation}
\end{lemma}
\begin{pf}
By the chain rule, (\ref{ET5}) follows by the following estimate:
%
\begin{equation}
\sup_{t\in[0,1]}\sup_{x\in\mathbb{R}^d}\bigl|\nabla^m
X'_t(x)\bigr|\leq C_{m,b},\label{ES3}
\end{equation}
which has been proved in estimating (\ref{ET6}).
\end{pf}
%
%
\begin{lemma}\label{Le5}
Let\vspace*{1pt} $J'_t(x):=\nabla X'_t(x)$ and $K'_t(x)$ be the inverse matrix of~$J'_t(x)$.
Let $f=(f_{kl})\in C^\infty_b(\mathbb{R}^d)$ be an $\mathbb
{R}^m\times\mathbb{R}^m$ valued function.\vspace*{1pt}
Then for any \mbox{$j=1,\ldots,d$} and $k,l=1,\ldots, m$, we have
the following formula:
%
\begin{equation}
\mathcal{P}'_t(\partial_j
f_{kl}) (x)=\operatorname{div}Q^{\cdot
j}_{kl}(t,x;
f)-G^j_{kl}(t,x; f),\label{ET9}
\end{equation}
where
%
\begin{eqnarray}
Q^{ij}_{kl}(t,x;f)&:=&\mathbb{E}\bigl(f_{kl}
\bigl(X'_t(x)\bigr) \bigl(K'_t(x)
\bigr)_{ij}\bigr),\label
{EE5}
\\
G^j_{kl}(t,x;f)&:=&\mathbb{E}\bigl(f_{kl}
\bigl(X'_t(x)\bigr)\operatorname {div}
\bigl(K'_t\bigr)_{\cdot j}(x)\bigr).\label{EE6}
\end{eqnarray}
Moreover, for any $m\in\{0\}\cup\mathbb{N}$, we have
%
\begin{equation}
\sup_{t\in[0,1]}\sup_{x\in\mathbb{R}^d}\bigl|\nabla^mK'_t(x)\bigr|
\leq \widetilde C_{m,b},\label{ES2}
\end{equation}
where $\widetilde C_{m,b}\geq1$.
\end{lemma}
\begin{pf}
Noticing that
\[
\nabla\bigl(f\bigl(X'_t(x)\bigr)\bigr)=(\nabla f)
\bigl(X'_t(x)\bigr)\nabla X'_t(x)=(
\nabla f) \bigl(X'_t(x)\bigr)J'_t(x),
\]
we have
\[
(\nabla f) \bigl(X'_t(x)\bigr)=\nabla\bigl(f
\bigl(X'_t(x)\bigr)\bigr) K'_t(x)=
\operatorname {div}\bigl(f\bigl(X'_t\bigr)
K'_t\bigr) (x)-f\bigl(X'_t(x)
\bigr)\operatorname{div}K'_t(x),
\]
which in turn gives (\ref{ET9}) by taking expectations.
As for (\ref{ES2}), it follows by equation
\[
K'_t(x)=I-\int^t_0K'_s(x)
\cdot\nabla b\bigl(X'_s(x)\bigr)\,\mathrm{d}s
\]
and estimate (\ref{ES3}).
\end{pf}

Below, let $\mathscr{C}:=\{\tau_1, \tau_2,\ldots,\tau_n,\ldots\}$
and $\mathscr{G}:=\{\xi_1,\xi_2,\ldots,\xi_n,\ldots\}$ be two
independent families of i.i.d. random variables
in $\mathbb{R}^+$ and $\mathbb{R}^d$, respectively,
which are also independent of $(W_t, S'_t)_{t\geq0}$. We assume that
$\tau_1$ obeys
the exponential distribution of parameter
\[
\lambda:=\nu_S\bigl([1,\infty)\bigr)
\]
and $\xi_1$ has the distributional density
\[
\frac{1}{\nu_S([1,\infty))}\int^\infty_1(2\pi
s)^{-d/2}\mathrm {e}^{-|x|^2/2s}\nu_S(\mathrm{d}s).
\]
Set $\tau_0:=0$ and $\xi_0:=0$, and define
\[
N_t:=\max\{n\dvtx  \tau_0+\tau_1+\cdots+
\tau_n\leq t\}=\sum_{n=0}^{\infty}
1_{\{\tau_0+\cdots+\tau_n\leq t\}}
\]
and
\[
H_t:=\xi_0+\xi_1+\cdots+
\xi_{N_t}=\sum_{j=0}^{N_t}
\xi_j.
\]
Then $H_t$ is a compound Poisson process with L\'evy measure
\[
\nu_H(\Gamma)=\int^\infty_1(2\pi
s)^{-d/2} \biggl(\int_\Gamma \mathrm{e}^{-|y|^2/2s}
\,\mathrm{d}y \biggr)\nu_S(\mathrm{d}s).
\]
Moreover, it is easy to see that $H_t$ is independent of $W_{S'_t}$, and
%
\begin{equation}
(AW_{S_t})_{t\geq0}\stackrel{(d)} {=}(AW_{S'_t}+AH_t)_{t\geq0}.\label{ET2}
\end{equation}

Let $\hbar_t$ be a c\`adl\`ag purely discontinuous $\mathbb
{R}^d$-valued function with finite many jumps and $\hbar_0=0$.
Let $X^\hbar_t(x)$ solve the following SDE:
\[
X^\hbar_t(x)=x+\int^t_0b
\bigl(X^\hbar_s(x)\bigr)\,\mathrm{d}s+AW_{S'_t}+
\hbar_t.
\]
Let $n$ be the jump number of $\hbar$ before time $t$. Let
$0=t_0<t_1<t_2<\cdots<t_n<t$ be the jump time of $\hbar$.
By the Markovian property of $X^\hbar_t(x)$, we have the following\vadjust{\goodbreak} formula:
\begin{eqnarray*}
&& \mathbb{E}f\bigl(X^\hbar_t(x)\bigr)
\\
&&\qquad = \int
_{\mathbb{R}^d} \biggl(\int_{\mathbb
{R}^d} \biggl(\int
_{\mathbb{R}^d}\cdots \biggl(\int_{\mathbb{R}^d}
p'_{t_1}(x,y_1)p'_{t_2-t_1}(y_1+
\Delta \hbar_{t_1},y_2)\,\mathrm{d}y_1 \biggr)
\\
&&\hspace*{129.5pt}\cdots p'_{t_n-t_{n-1}}(y_{n-1}+\Delta
\hbar_{t_{n-1}},y_n) \,\mathrm{d}y_{n-1} \biggr)
\\
&&\hspace*{184pt}{}\times p'_{t-t_n}(y_n+\Delta\hbar_{t_n},z)
\,\mathrm{d}y_n \biggr)f(z)\,\mathrm{d}z
\\
&&\qquad =\mathcal{P}'_{t_1}\cdots
\vartheta_{\Delta\hbar_{t_{n-1}}} \mathcal{P}'_{t_n-t_{n-1}}
\vartheta_{\Delta\hbar_{t_{n}}}\mathcal {P}'_{t-t_{n}}f(x),
\end{eqnarray*}
where
\[
\vartheta_y g(x):=g(x+y).
\]
Now, by (\ref{ET2}) we have
\[
X_t(x)\stackrel{(d)} {=}X^\hbar_t(x)|_{\hbar=AH_\cdot}
\]
and so,
\begin{eqnarray*}
\mathcal{P}_tf(x)&=&\mathbb{E}f\bigl(X_t(x)\bigr)=
\mathbb{E} \bigl(\mathbb {E}f\bigl(X^\hbar_t(x)\bigr)
|_{\hbar=AH_\cdot} \bigr)
\\
&=&\sum_{n=0}^\infty\mathbb{E} \bigl(
\mathcal{P}'_{\tau_1}\cdots \vartheta_{A\xi_{n-1}}
\mathcal{P}'_{\tau_n}\vartheta_{A\xi_{n}}
\mathcal{P}'_{t-(\tau
_0+\tau_1+\cdots+\tau_n)}f(x); N_t=n \bigr).
\end{eqnarray*}
In view of
\[
\{N_t=n\}=\{\tau_0+\cdots+\tau_n\leq t<
\tau_0+\cdots+\tau_{n+1}\}
\]
and that $\mathscr{C}$ is independent of $\mathscr{G}$, we further have
\begin{eqnarray}\label{For}
\mathcal{P}_tf(x)\nonumber
&=&\sum_{n=1}^\infty
\biggl\{\int_{t_1+\cdots
+t_n<t<t_1+\cdots+t_{n+1}} \lambda^{n+1}
\mathrm{e}^{-\lambda(t_1+\cdots
+t_n+t_{n+1})}\hspace*{-25pt}\nonumber
\\
&&\hspace*{27pt}{} \times\mathbb{E} \bigl(\mathcal{P}'_{t_1}\cdots
\vartheta _{A\xi_{n-1}} \mathcal{P}'_{t_n}
\vartheta_{A\xi_{n}}\mathcal{P}'_{t-(t_1+\cdots
+t_n)}f(x) \bigr)
\,\mathrm{d}t_1\cdots\mathrm{d}t_{n+1} \biggr\}\hspace*{-25pt}
\nonumber
\\
&&{}+\mathcal{P}'_tf(x)P(N_t=0)\hspace*{-25pt}
\\
&=&\sum_{n=1}^\infty \biggl\{
\lambda^{n} \mathrm{e}^{-\lambda t}\int_{t_1+\cdots+t_n<t}
\mathbb{E}I^{A\bolds{\xi}}_f(t_1,\ldots,
t_n, t,x) \,\mathrm{d}t_1\cdots\mathrm{d}t_{n}
\biggr\}\hspace*{-25pt}\nonumber
\\
&&{}+\mathcal{P}'_tf(x)\mathrm{e}^{-\lambda t},\hspace*{-25pt}\nonumber
\end{eqnarray}
where $\bolds{\xi}:=(\xi_1,\ldots,\xi_n)$, and
\[
I^{\mathbf{y}}_f(t_1,\ldots, t_n,
t,x):=\mathcal{P}'_{t_1}\cdots \vartheta_{y_{n-1}}
\mathcal{P}'_{t_n}\vartheta_{y_n}
\mathcal{P}'_{t-(t_1+\cdots+t_n)}f(x)
\]
with ${\mathbf{y}}:=(y_1,\ldots,y_n)$.

Now we can complete the proof of Theorem~\ref{Th2}.

\begin{pf*}{Proof of Theorem~\ref{Th2}}
We first establish the same gradient
estimate as in~(\ref{ET7}).

If we let $t_{n+1}:=t-(t_1+\cdots+t_n)>0$, then there is at least one
$j\in\{1,2,\ldots,n+1\}$ such that
%
\begin{equation}
t_j\geq\frac{t}{n+1}.\label{Eq3}
\end{equation}
Thus, we have
\begin{eqnarray*}
\bigl|\nabla_xI^{\mathbf{y}}_f(t_1,
\ldots,t_n,t,x)\bigr|&\stackrel{\fontsize{8.36}{10}{\selectfont{(\ref{ET5})}}} {\leq}&
C_{1,b}^{j-1} \bigl\llVert \nabla_x
\mathcal{P}'_{t_j}\cdots\vartheta_{y_{n-1}}
\mathcal{P}'_{t_n}\vartheta_{y_{n}}
\mathcal{P}'_{t_{n+1}}f\bigr\rrVert _\infty
\\
&\stackrel{\fontsize{8.36}{10}{\selectfont{(\ref{ET7})}}}{\leq}& C C_{1,b}^{j-1}t_j^{-{\gamma
_{1,0}}}
\bigl\llVert \mathcal{P}'_{t_{j+1}}\cdots\vartheta_{y_{n-1}}
\mathcal{P}'_{t_n}\vartheta_{y_{n}}
\mathcal{P}'_{t_{n+1}}f\bigr\rrVert _\infty
\\
&\leq& C C_{1,b}^n\bigl(t/(n+1)\bigr)^{-\gamma_{1,0}}\|f\|_\infty.
\end{eqnarray*}
Here and below, the various constant $C$ is independent of $t$ and $n$.
Hence, by (\ref{For}) we have
%
\begin{eqnarray}\label{EE4}
\qquad \bigl|\nabla\mathcal{P}_tf(x)\bigr|&\leq& C\|f\|_\infty
t^{-\gamma _{1,0}}\mathrm{e}^{-\lambda t}
\nonumber
\\
&&{} \times \Biggl(1+\sum_{n=1}^\infty
\lambda^{n} C_{1,b}^n(n+1)^{\gamma_{1,0}} \int
_{t_1+\cdots+t_n<t}\,\mathrm{d}t_1\cdots \mathrm{d}t_{n}
\Biggr)
\nonumber\\[-8pt]\\[-8pt]
&=&C\|f\|_\infty t^{-\gamma_{1,0}}\mathrm{e}^{-\lambda t} \Biggl(\sum
_{n=0}^\infty\lambda^{n}
C_{1,b}^n(n+1)^{\gamma
_{1,0}}\frac{t^n}{n!} \Biggr)\nonumber
\\
&\leq& C\|f\|_\infty t^{-\gamma_{1,0}}.\nonumber
\end{eqnarray}
Thus, we obtain (\ref{ET7}) with $k=1$ and $m=0$.

For $k,l=1,\ldots,d$, set $F^{(0)}_{kl}(x):=1_{k=l}f(x)$ and
$R^{(0)}_{l}(x):=0$. Let us recursively define for $m=0,1,\ldots,n$,
\begin{eqnarray*}
F^{(m+1)}_{kl}(x)&:=&\sum_{i=1}^d
Q^{ki}_{il}\bigl(t_{n+1-m},x;\vartheta
_{y_{n+1-m}}F^{(m)}\bigr),
\\
R^{(m+1)}_{l}(x)&:=&\sum_{i=1}^dG^{i}_{il}
\bigl(t_{n+1-m},x;\vartheta _{y_{n+1-m}}F^{(m)}\bigr),
\end{eqnarray*}
where $y_{n+1}:=0$, $Q^{ki}_{il}$ and $G^{i}_{il}$ are defined by (\ref
{EE5}) and (\ref{EE6}).
From these definitions and by (\ref{ES2}), it is easy to see that
\begin{eqnarray*}
\bigl\|F^{(m+1)}_{kl}\bigr\|_\infty&\leq& d
\bigl\|F^{(m)}_{kl}\bigr\|_\infty\mathbb{E}\bigl\|
K'_t(x)\bigr\|\leq\widetilde C_{0,b}
\bigl\|F^{(m)}_{kl}\bigr\|_\infty
\\
&\leq&\widetilde C_{0,b}^{m+1}\bigl\|F^{(0)}_{kl}
\bigr\|_\infty\leq\widetilde C_{0,b}^{m+1}\|f
\|_\infty
\end{eqnarray*}
and
\[
\bigl\|R^{(m+1)}_{l}\bigr\|_\infty\leq\bigl\|F^{(m)}_{kl}
\bigr\|_\infty\mathbb {E}\bigl|\operatorname{div}K'_t(x)\bigr|
\leq\widetilde C_{0,b}^m\widetilde C_{1,b}\|f\|
_\infty.
\]\eject\noindent
By repeatedly using Lemma~\ref{Le5}, we have
\begin{eqnarray*}
&& \bigl|I^{\mathbf{y}}_{\partial_l f}(t_1,\ldots,t_n,t,x)\bigr|
\\
&&\quad =
\Biggl|\mathcal{P}'_{t_1}\cdots\vartheta_{y_{j-1}}
\mathcal {P}'_{t_j}\operatorname{div}F^{(n+1-j)}_{\cdot l}(x)
-\sum_{m=1}^{n+1-j}
\mathcal{P}'_{t_1}\cdots\vartheta _{y_{n+1-m}}
\mathcal{P}'_{t_{n+1-m}}R^{(m)}_{l}(x) \Biggr|
\\
&&\hspace*{6pt}\stackrel{\fontsize{8.36}{10}{\selectfont{(\ref{ET7})}}} {\leq} C t_j^{-\gamma_{0,1}}\bigl\|
F^{(n+1-j)}_{\cdot l}\bigr\|_\infty+\sum
_{m=1}^{n+1-j}\bigl\|R^{(m)}_{l}\bigr\|
_\infty
\\
&&\hspace*{6pt}\stackrel{\fontsize{8.36}{10}{\selectfont{(\ref{Eq3})}}} {\leq}C\bigl(t/(n+1)\bigr)^{-\gamma_{0,1}}\widetilde
C_{0,b}^{n}\|f\|_\infty+C\widetilde
C_{0,b}^{n}\|f\|_\infty.
\end{eqnarray*}
As in estimating (\ref{EE4}), we obtain (\ref{ET7}) with $k=0$ and
$m=1$. For the general $m$ and $k$,
the gradient estimate (\ref{ET7}) follows by similar calculations and
the induction.

Lastly, by estimate (\ref{ET7}) and Sobolev's embedding theorem (see
\cite{Nu}, \mbox{pp. 102--103}),
one has that for each $t>0$,
\[
(x,y)\mapsto p_t(x,y)\in C^\infty_b\bigl(
\mathbb{R}^d\times\mathbb{R}^d\bigr).
\]
The smoothness of $p_t(x,y)$ with respect to the time variable $t$ follows
by equation~(\ref{ET8}) and the standard bootstrap argument. As for
equation (\ref{Eq2}), it follows~by
\[
\frac{\mathrm{d}\mathcal{P}_t f(x)}{\mathrm{d}t}=\mathcal {L}_A\mathcal{P}_t
f(x)+b(x)\cdot\nabla_x\mathcal{P}_tf(x),
\]
where $f\in C^\infty_b(\mathbb{R}^d)$.
\end{pf*}

\section*{Acknowledgements}
The author is very grateful to Professors Hua Chen, Zhen-Qing Chen and
Feng-Yu Wang for their quite useful conversations.
The referees' very useful suggestions are also deeply acknowledged.




\printaddresses

\end{document}